

\documentstyle{amsppt}                 
\NoBlackBoxes

 \def\a{\kern+.6ex\lower.42ex\hbox{$\scriptstyle \iota$}\kern-1.20ex a}
 \def\e{\kern+.5ex\lower.42ex\hbox{$\scriptstyle \iota$}\kern-1.10ex e}
 \define\Ce{\Bbb C}                         
 \define\Ne{\Bbb N}                         
 \define\ord{\operatorname{ord}}            
 \define\be{\bar b}                         
 \define\series#1{\Ce\bigl[ [ #1 ]\bigr]}   

\topmatter
 \title On the approximate roots of polynomials \endtitle
 \author by Janusz Gwo\'zdziewicz and Arkadiusz P\l{}oski \endauthor
 \address
 \newline
   Department of Mathematics,\newline
   Technical University,\newline
   Al.~1000\,LPP\,7, 25--314~Kielce,\newline
   Poland
 \endaddress
 \email mat-jg\@srv1.tu.kielce.pl \endemail

 \abstract
   We give a simplified approach to the Abhyankar--Moh theory of
   approximate roots. Our considerations are based on properties
   of the intersection multiplicity of local curves.
  \endabstract
\endtopmatter

\document
\head 		1. The main results
\endhead        
For any power series $f$, $g\in\series{x,y}$ we define the
intersection number  $(f,g)_0$ by
$(f,g)_0=\dim_{\Ce}\series{x,y}/{(f,g)}$.
  Suppose that $f=f(x,y)$ is an irreducible power series and let
$n=(f,x)_0=\ord f(0,y) < {+}\infty$.
Then there exists a power series $y(t)\in\series{t}$,
$\ord y(t)>0$ such that $f(t^n,y(t))=0$.
We have $(f,g)_0=\ord g(t^n,y(t))$ for any $g=g(x,y)\in\series{x,y}$.
The mapping $g\mapsto (f,g)_0$ induces a valuation $v_f$ of the
ring $\series{x,y}/(f)$. Let $\Gamma(f)$ be the semigroup of
$v_f$ i.e.{}
$\Gamma (f) = \{\,(f,g)_0\in\Ne : g\not\equiv 0 \mod (f)\,\}$.

According to the well known structure theorem for the semigroup
$\Gamma (f)$ (\cite{1}, \cite{2}, \cite{10} and Sect\. ~3 of
this paper) there is a unique sequence of positive integers
$\be_0$, $\be_1$,\dots , $\be_h$ such that

 \item{(i)}
     $\be_0 = (f,x)_0$,

 \item{(ii)}
     $\be_k=\min(\Gamma (f)\setminus (\Ne\be_0+\dots +\Ne\be_{k-1}))$
     \quad for $k=1,\dots ,h$.

 \item{(iii)}
     $\Gamma (f)=\Ne\be_0+\dots +\Ne\be_h$ \quad i.e.~
     $\Gamma (f)$ is generated by
     $\be_0, \be_1,\dots , \be_h$

If the conditions (i), (ii), (iii) are satisfied, we write
$\Gamma (f)=\langle \be_0,\dots ,\be_h \rangle$.
Also define
$B_k = \gcd (\be_0,\dots ,\be_k)$ \, for $k=0, 1,\dots, h$
and $n_k=B_{k-1}/B_{k}$ \, for $k=1,\dots, h$. We have

\item{(iv)}
    $n_k>1$ and the sequence $B_{k-1}\be_k$ is strictly
    increasing for $k\geq 1$.

Let $\Cal O$ be an integral domain of characteristic zero.
Let $g\in\Cal O[y]$ be a monic polynomial and let $d$ be a
positive integer such that $d$ divides $\deg g$. According to
Abhyankar and Moh \cite{2} the approximate $d^{{\text{th}}}$ root
of $g$ denoted $\root{d}\of{g}$ is defined to be the unique
monic polynomial satisfying
$\deg (g-(\root{d}\of{g})^d) < \deg g - \deg \root{d}\of{g}$.
Obviously $\deg \root{d}\of{g} = \deg g/d$.

Let $1\leq k\leq h$.

\proclaim{Theorem 1.1}
%
Let $g=g(x,y)\in\series{x}[y]$ be a monic polynomial,
$\deg_y g=n/B_k$. Suppose that $(f,g)_0>n_k\be_k$.
Then $(f,\root{n_k}\of{g})_0 = \be_k$.
\endproclaim

The proof of (1.1) is given in Sect\.~4 of this paper.
We shall follow the methods developed by Abhyankar and Moh
in the fundamental paper ~\cite{2} and simplified by Abhyankar
in ~\cite{1}.

Let $1\leq k\leq h+1$.
\proclaim{Theorem 1.2}
%
Let $\phi=\phi(x,y)\in\series{x,y}$ be such that
$(\phi,x)_0 = n/B_{k-1}$.
Then $(f,\phi)_0\leq\be_k$ (we put $\be_{h+1}={+}\infty$).
If additionaly, $(f,\phi)_0 > n_{k-1}\be_{k-1}$ then $\phi$ is
irreducible and
$\Gamma (\phi) = \langle
    \frac{\be_0}{B_{k-1}},\dots ,\frac{\be_{k-1}}{B_{k-1}}
                 \rangle$
\endproclaim

The proof of (1.2) is given in Sect\. ~5 of this paper.
Note that the second part of (1.2) is a
generalization of Abhyankar's irreducibility criterion
\cite{1, Chapter V}.

\proclaim{Lemma~1.3}
%
Let $g\in\Cal O[y]$ be a monic polynomial of degree $n>0$ and
let $d, e>0$ be integers such that $de$ divides $n$.
Then $\root{e}\of{\root{d}\of{g}} = \root{ed}\of{g}$.
\endproclaim

\demo{Proof of (1.3)}
Let $h=\root{e}\of{\root{d}\of{g}}$, so $\deg h=\frac{n}{de}$.
We have
$\root{d}\of{g} = (\root{e}\of{\root{d}\of{g}})^e + R$,\quad
$\deg R<\frac{n}{d}-\frac{n}{de}$
and
$g=(\root{d}\of{g})^d+S$, $\deg S<n-\frac nd$.
Then we get $g=(h^e+R)^d+S = h^{ed}+T$ where
$T=\sum_{i=1}^d \binom{d}{i} (h^e)^{d-i}R^i+S$.
Since
$ \deg (\binom{d}{i} (h^e)^{d-i}R^i) <
  (ed-ei)\frac{n}{de}+i(\frac{n}{d}-\frac{n}{de}) \leq
  n-\frac{n}{de}
$
then $\deg T < n-\frac{n}{de}$.
This shows that $h=\root{ed}\of{g}$.
\enddemo

Now, we can prove the Abhyankar--Moh theorem.
\proclaim{Theorem~1.4 \cite{2}}
%
Let $f=f(x,y)\in\series{x}[y]$ be an irreducible distinguished
polynomial of degree $n>1$ with
$\Gamma (f) = \langle \be_0,\be_1,\dots ,\be_h \rangle$ and
$\be_0=(f,x)_0=n$.
Let $1\leq k\leq h+1$.
Then we have
\item{(1.4.1)}     $(f,\root{B_{k-1}}\of{f})_0 = \be_k$

\item{(1.4.2)}
    $\root{B_{k-1}}\of{f}$ is an irreducible distinguished polynomial
    of degree $\frac{n}{B_{k-1}}$ such that
    $\Gamma (\root{B_{k-1}}\of{f}) =
     \langle
       \frac{\be_0}{B_{k-1}},\dots ,\frac{\be_{k-1}}{B_{k-1}}
     \rangle$.
\endproclaim

\demo{Proof of (1.4)}
First we check (1.4.1) using (1.1) and induction on $k$.
If $k=h+1$ then $B_{k-1}=B_h$, $\be_k=\be_{h+1}={+}\infty$ and
(1.4.1) is obvious.
Let $k\leq h$ and suppose $(f,\root{B_k}\of{f})_0=\be_{k+1}$
holds true.

The polynomial $\root{B_k}\of{f}$ is of degree $n/B_k$ and
$(f,\root{B_{k-1}}\of{f})_0 > n_k\be_k$\quad
($\be_{k+1} > n_k\be_k$ by ~(iv))
so we can apply ~(1.1) to $g=\root{B_k}\of{f}$.
By ~(1.1) we get $(f,\root{n_k}\of{g})_0 = \be_k$ and by ~(1.3) we have
$\root{n_k}\of{g} = \root{n_k}\of{\root{B_k}\of{f}}=\root{B_{k-1}}\of{f}$,
consequently $(f,\root{B_{k-1}}\of{f})_0 = \be_k$ and ~(1.4.1) is proved.

In order to prove (1.4.2) put $\phi = \root{B_{k-1}}\of{f}$. Thus
$\phi$ is a distinguished polynomial of degree $n/B_{k-1}$, hence
$(\phi,x)_0 = n/B_{k-1}$. On the other hand
$(f,\phi)_0=\be_k>n_{k-1}\be_{k-1}$ by (1.4.1).
According to theorem~1.2 $\phi = \root{B_{k-1}}\of{f}$ is irreducible
and
$\Gamma (\phi) =
        \langle
            \frac{\be_0}{B_{k-1}},\dots ,\frac{\be_{k-1}}{B_{k-1}}
        \rangle$.
\enddemo

\remark{Note}
Theorems (1.1) and (1.2) could be formulated and proved for meromorphic
curves. However we shall see in the next section that the algebroid
case considered by us is sufficient to get important applications
of approximate roots.
\endremark

\head 		2. The Abhyankar--Moh inequality
\endhead        
We give here a geometrical version of the Abhyankar--Moh inequality
which is the basic tool for proving the Embedding Theorem
\cite{3}. Let $\bold C\subset\Bbb P^2$ be an irreducible projective
plane curve of degree $n>1$ and let $\bold O\in\bold C$ be its
singular point. We assume that $\bold C$ is analytically
irreducible at $\bold O$ i.e\. the analytic germ
$(\bold C,\bold O)$ is irreducible,
and let $\bold L$ be the unique tangent to
$\bold C$ at $\bold O$. Let $\Gamma(\bold C,\bold O)$ be the
semigroup of the branch of $\bold C$ passing through $\bold O$
and let
$\Gamma(\bold C,\bold O) = \langle \be_0,\be_1, \dots ,\be_h \rangle$
with $\be_0 = (\bold C \centerdot \bold L)_{\bold O}$.

\proclaim{Theorem~2.1 \cite{3}}
%
Suppose that $\bold C\cap\bold L=\{\bold O\}$, i.e. $\be_0=n$.
Then we have $B_{h-1}\be_h<n^2$.
\endproclaim
\demo{Proof}
Choose the line at infinity no passing through $\bold O$
and let $x$, $y$
be an affine system of coordinates centered at $\bold O$ such
that $\bold L$ has an equation $x=0$. Let $f(x,y)\in\Ce [x,y]$
be the irreducible equation of $\bold C$ in coordinates $x$, $y$.
It is easy to see that $f$ is a distinguished, irreducible in
$\series{x,y}$ polynomial and
$\Gamma (f) = \langle \be_0,\be_1, \dots ,\be_h \rangle$.

We have $\deg f=n$ and consequently
$\deg\root{B_{h-1}}\of{f}=\frac{n}{B_{h-1}}$, thus (1.4.1) and
Bezout's theorem imply
$\be_h=(f,\root{B_{h-1}}\of{f})_0\leq n\frac{n}{B_{h-1}}$ that
is $B_{h-1}\be_h\leq n^2$.
In fact $B_{h-1}\be_h<n^2$ because the equality
$B_{h-1}\be_h=n^2$ implies $\be_h=n\frac{n}{B_{h-1}}$ which
contradicts the relation $\be_h\not\equiv 0\mod B_{h-1}$.
\enddemo

The inequality (2.1) has an application to the polar
curves \cite{4}.

\proclaim{Theorem~2.2 (with the assumptions as above)}
%
Let $(\bold D,\bold O)$ be an irreducible component of the local
polar of $\bold C$ with respect to $\bold L$.

Then
$(\bold C\centerdot\bold D)_{\bold O} <
 (\bold C\centerdot\bold L)_{\bold O}
 (\bold D\centerdot\bold L)_{\bold O}$.
\endproclaim
\demo{Proof}
In the coordinates $x$, $y$ introduced in the proof of (2.1) the
local polar is given by equation
$\frac{\partial f}{\partial y}=0$ and its irreducible component
is given by $g=0$ where $g$ is irreducible (in $\series{x,y}$)
divisor of  $\frac{\partial f}{\partial y}$. By the Merle
formula for polar invariants \cite{6}, \cite{4}, \cite{5}
and theorem~2.1 we get
$$
\frac{(\bold C\centerdot\bold D)_{\bold O}}
     {(\bold D\centerdot\bold L)_{\bold O}} =
 \frac{(f,g)_0}{(g,x)_0} = \frac{B_{k-1}}{\be_0}\be_k \leq
 \frac{B_{h-1}}{\be_0}\be_h < n = (f,x)_0 =
 (\bold C\centerdot\bold L)_{\bold O}
$$
and the theorem follows.
\enddemo

\head         3. Characteristic, semigroup of an analytic curve
                                   and  the Noether formula
\endhead      
In this section we recall some well-known notions of the theory
of analytic curves. Our main reference is \cite{10}.
Let $f=f(x,y)$ be an irreducible power series $y$--regular of
order $n=\ord f(0,y)>1$. There exists a power series
$y(t)\in\series{t}$, $\ord y(t)>0$ such that $f(t^n,y(t))=0$.
Moreover every solution of the equation $f(t^n,y)=0$ is of the
form $y(\epsilon t)$ for some $\epsilon$ such that
$\epsilon^n=1$. Let $y(t)=\sum a_jt^j$.
We put $S(f) = \{\,j\in\Ne : a_j\neq 0\,\}$.  Note that $S(f)$
depends only on $f$.

The characteristic $b_0, b_1, \dots, b_h$ of $f$ is the unique
sequence of positive integers satisfying the conditions:
 \item{(i)}
         $b_0=n$
 \item{(ii)}
         $b_{k+1}=\min\{\, j\in S(f) : \gcd (b_0,\dots ,b_k,j)<
         \gcd (b_0,\dots ,b_k)\,\}$
 \item{(iii)}
         $\gcd (b_0,\dots ,b_h)=1$

In the sequel we put
$B_k=\gcd (b_0,\dots ,b_k)$\quad for $k=0$, $1$,\dots, $h$

and
$\be_k = b_k + \frac{1}{B_{k-1}}\sum_{i=1}^{k-1}(B_{i-1}-B_i)b_i$\quad
for $k=1$,\dots, $h$.

We assume that the sum of an empty family is equal to zero. Thus
we have $\be_1=b_1$. We put $\be_0=b_0$. One checks easily that
$\gcd(\be_0,\dots ,\be_k)=B_k$ for $k=0$,\dots, $h$ and
$  b_k = \be_k - \sum_{i=1}^{k-1}(\frac{B_{i-1}}{B_i}-1)\be_i$
\quad for $k=0$, $1$,\dots, $h$.
Therefore the sequence ~$\be_0$,\dots, $\be_k$ determines
sequences $B_0$,\dots, $B_k$ and $b_0$,\dots, $b_k$ for any
$k=0$,\dots ,$h$.
We have
$\be_{k+1}-\frac{B_{k-1}}{B_k}\be_k=b_{k+1}-b_k$
for $k=1$,\dots, $h-1$, which shows that the sequence
{}~$B_{k-1}\be_k$ is increasing.

For any $k$, $1\leq k\leq h$ we set
$S_k=\{\,j\in S(f):j<b_{k+1}\,\}$
We put $b_{h+1}=\be_{h+1}={+}\infty$, so $S_h=S(f)$.  Let
$y_k(t)=\sum\limits_{j\in S_k}a_jt^j$. There exists an
irreducible, monic polynomial $f_k=f_k(x,y)\in\series{x}[y]$
such that $f_k(t^n,y_k(t))=0$.

\proclaim{Lemma~3.1 \cite{10}}
%
 $\deg_yf_k(x,y)=n/B_k$,\quad $(f,f_k)_0=\be_{k+1}$
\endproclaim
\demo{Proof} \cite{10, p\. ~15} \enddemo

Note that $f_{h+1}$ is the distinguished polynomial associated
with $f$.

\proclaim{Proposition~3.2 \cite{10}}
%
If $\psi(x,y)\in\series{x}[y]$, $\deg_y\psi(x,y)<n/B_k$ and
$\psi\not\equiv 0 \mod (f)$, then
$(f,\psi)_0\in\Ne\be_0+\dots+\Ne\be_k$
\endproclaim
\demo{Proof} \cite{10, p\. ~16} \enddemo

Note that from (3.1) and (3.2) we get
$\Gamma(f)=\Ne\be_0+\dots+\Ne\be_h$. Now, let $g=g(x,y)$ be an
irreducible power series $y$--regular of order $p=\ord
f(0,y)<{+}\infty$. Suppose that $f$, $g$ are coprime. Let
$z(t)\in\series{t}$, $\ord z(t)>0$ be such that $g(t^p,z(t))=0$.
We put
$$
o_f(g)=\max\{\,\ord(y(\epsilon x^{1/n})-z(\nu x^{1/p})) :
               \epsilon^n=1, \nu^p=1\,\}.
$$
It is easy to check, that
$$\align
o_f(g) = &\max\{\,\ord(y(x^{1/n})-z(\nu x^{1/p})) : \nu^p=1\,\}     \\
    {} = &\max\{\,\ord(y(\epsilon x^{1/n})-z(x^{1/p})) : \epsilon^n=1\,\}
\endalign
$$
In particular $o_f(g)=o_g(f)$.

The classical computation leads to the following formula due to
Max Noether:

\proclaim{Proposition~3.3 \cite{6}, \cite{5}}
%
Suppose that $f$, $g$ are irreducible, $y$--regular power series,
$f$ of characteristic $(b_0,b_1,\dots, b_h)$ and let $k$ be the
smallest strictly positive integer such that
$o_f(g)\leq \frac{b_k}{b_0}$ \,
($\frac{b_{h+1}}{b_0}={+}\infty$).
Then
$$
  \frac{(f,g)_0}{(g,x)_0} =
  \sum_{i=1}^{k-1}(B_{i-1}-B_i)\frac{b_i}{b_0}\,+B_{k-1}o_f(g)
$$
\endproclaim

\remark{Remark \cite{9}}
The Noether formula is really symetric.
Let $(c_0,c_1,\dots,c_m)$ be the characteristic of $g$. Then
$k\leq m$,\quad $\frac{c_i}{c_0}=\frac{b_i}{b_0}$ for $i=1$,\dots,
$k-1$  and $o_f(g)\leq\frac{c_k}{c_0}$. If
$C_i=\gcd(c_0,\dots ,c_i)$ then the formula can be rewritten in
the following form
$$
  \frac{(f,g)_0}{(f,x)_0} =
  \sum_{i=1}^{k-1}(C_{i-1}-C_i)\frac{c_i}{c_0}\,+C_{k-1}o_g(f).
$$
\endremark

Using (3.3) we check easily
\proclaim{Lemma~ 3.4 }
%
Let $l>0$ be an integer. Then $o_f(g)\leq \frac{b_l}{b_0}$
iff
$\frac{(f,g)_0}{(g,x)_0}\leq B_{l-1}\frac{\be_l}{\be_0}$.
Moreover $o_f(g)=\frac{b_l}{b_0}$
is equivalent to
$\frac{(f,g)_0}{(g,x)_0}=B_{l-1}\frac{\be_l}{\be_0}$.
\endproclaim

To end with, let us note
\proclaim{Corollary 3.5 to theorem 1.2}
%
If $f=f(x,y)\in\series{x}[y]$ is an irreducible distinguished
polynomial, then $o_f(\root{B_{k-1}}\of{f})=\frac{b_k}{b_0}$
\quad for $k=1$, \dots, $h$.
\endproclaim
\demo{Proof}
We have
$(\root{B_{k-1}}\of{f},x)_0=\deg_y\root{B_{k-1}}\of{f}=\frac{n}{B_{k-1}}$
and $(f,\root{B_{k-1}}\of{f})_0=\be_k$ by (1.4.1). Hence
$\dfrac{(f,\root{B_{k-1}}\of{f})_0}{(\root{B_{k-1}}\of{f},x)_0}=
 B_{k-1}\dfrac{\be_k}{\be_0}$.
The power series $\root{B_{k-1}}\of{f}$ is irreducible
by~(1.4.2), thus we get by (3.4) \quad
$o_f(\root{B_{k-1}}\of{f})=\frac{b_k}{b_0}$.
\enddemo

\head         4. Proof of theorem~1.1
\endhead      
Let $g\in\Cal O[y]$ be a monic polynomial with coefficients in
the integral domain $\Cal O$ of charecteristic zero and let $d$
be a positive divisor of $\deg g$. Given any monic polynomial
$h\in\Cal O[y]$ of degree $\deg g/d$ we have $h$--adic expansion
of $g$, namely
$$
   g=h^d+a_1h^{d-1}+\dots+a_d,\quad
   a_i\in\Cal O[y],\quad
   \deg a_i<\deg h
$$
The polynomials $a_i$ are uniquely determined by $g$, $h$. The
Tschirnhausen operator $\tau_g(h)=h+\frac{1}{d}a_1$ changes $h$
to $\tau_g(h)$ which is again monic of degree $\deg g/d$.

\proclaim{Lemma~4.1 \cite{1}}
%
$\root{d}\of{g}=\tau_g(\tau_g\dots(\tau_g(h))\dots)$ with
$\tau_g$ repeated $deg g/d$ times.
\endproclaim
\demo{Proof \cite{1, p\. 16}} \enddemo
To prove (1.1) it suffices to check the following

\item {(*)}
   if $h(x,y)\in\series{x}[y]$ is a monic polynomial of degree
   $n/B_{k-1}$ such that $(f,h)_0=\be_k$, then
   $(f,\tau_g(h))_0=\be_k$.

Indeed, to get the relation $(f,\root{n_k}\of{g})_0=\be_k$ we
take $h=f_{k-1}$ (cf\. lemma~ 3.1) and apply Tschirnhausen
operator $\tau_g$ to $h$\quad $\deg g/n_k=n/B_{k-1}$~ times.

To prove (*) fix a monic polynomial $h(x,y)\in\series{x}[y]$
such that $\deg h=n/B_{k-1}$ and $(f,h)_0=\be_k$ and let us
consider the $h$--adic expansion of $g$:
$$
     g=h^{n_k}+a_1h^{n_k-1}+\dots+a_{n_k},\qquad
       \deg_ya_i<\deg_yh=n/B_{k-1}              \tag 1
$$
Let $I$ be the set of all $i\in\{\,1,\dots, n_k\,\}$ such that
$a_i\neq 0$. Therefore $(f,a_i)_0<{+}\infty$ for $i\in I$ and
by Proposition~3.2 we have
$(f,a_i)_0\in\Ne\be_0+\dots+\Ne\be_{k-1}$ hence
$(f,a_i)_0\equiv 0\mod B_{k-1}$ for $i\in I$.

We have
$$
  (f,a_ih^{n_k-i})_0\neq(f,a_jh^{n_k-j})_0 \qquad
   \text{for } i\neq j \in I                         \tag 2
 $$
Suppose that (2) is not true, So there exist $i$, $j\in I$
such that $i<j$ and $(f,a_ih^{n_k-i})_0=(f,a_jh^{n_k-j})_0$.
Therefore $(f,a_i)_0+(n_k-i)\be_k=(f,a_j)_0+(n_k-j)\be_k$ and
$(j-i)\be_k=(f,a_j)_0-(f,a_i)_0\equiv 0\mod B_{k-1}$.
The last relation implies                                
$(j-i)\frac{\be_k}{B_k}\equiv 0\mod n_k$ and consequently
$j-i\equiv 0 \mod n_k$ because $\be_k/B_k$ and $n_k$ are coprime.
We get a contradiction because $0<j-i<n_k$.

{}From (1) and (2) we have
$$
    (f,g-h^{n_k})_0=min_{i=1}^{n_k}(f,a_ih^{n_k-i})_0     \tag 3
$$
By assumption $(f,g)_0>n_k\be_k=(f,h^{n_k})_0$, so
$(f,g-h^{n_k})_0=n_k\be_k$ and (3) implies
$n_k\be_k\leq(f,a_ih^{n_k-i})_0=(f,a_i)_0+(n_k-i)\be_k$ for
$i=1$,\dots, $n_k$.
Therefore we get
$$
  (f,a_i)_0\geq i\be_k \qquad \text{for } i=1,\dots, n_k \tag 4
$$
Moreover we have
$$
 \text{if}\quad (f,a_i)_0=i\be_k,\quad 1\leq i\leq n_k \qquad
 \text{then}\quad i=n_k .                 \tag 5
$$
Indeed, from $(f,a_i)_0=i\be_k$ it follows that
$i\be_k\equiv 0\mod B_{k-1}$ and
$i\frac{\be_k}{B_k}\equiv 0\mod n_k$, so $i\equiv 0\mod n_k$
because $\be_k/B_k$, $n_k$ are coprime. Hence we get $i=n_k$.

By (5) we get (because $n_k>1$)
$$ (f,a_1)_0>b_k    \tag 6   $$
Therefore
$(f,\tau_g(h))_0=(f,h+\frac{1}{n_k}a_1)_0=(f,h)_0=\be_k$.

\head 		5. Proof of theorem~1.2
\endhead        
The proof of (1.2) is based on the following

\proclaim{Lemma~5.1}
%
Let $g=g(x,y)$ be an irreducible power series,
$p=(g,x)_0<{+}\infty$ and let $1< k\leq h+1$.
If $\frac{(f,g)_0}{(g,x)_0}>\frac{B_{k-2}\be_{k-1}}{b_0}$,
then $(g,x)_0\equiv 0\mod \frac{b_0}{B_{k-1}}$.
If, additionally $(g,x)_0=\frac{b_0}{B_{k-1}}$, then
$\Gamma(f)=
  \langle
  \frac{\be_1}{B_{k-1}},\dots,\frac{\be_{k-1}}{B_{k-1}}
  \rangle$.
\endproclaim
\demo{Proof of 5.1}
Let $(c_0,c_1,\dots,c_m)$, $c_0=p$ be the characteristic of $g$.
By lemma~3.4 we have $o_f(g)>\frac{b_{k-1}}{b_0}$, so there
exist Puiseux expansions determined by $f(x,y)=0$ and $g(x,y)=0$
respectively which coincide up to the `monomials' of degree
$\frac{b_{k-1}}{n}$. Therefore $k-1\leq m$ and
$\frac{b_1}{n}=\frac{c_1}{p}$,\dots,
$\frac{b_{k-1}}{n}=\frac{c_{k-1}}{p}$.
There exist integers $a_0$,\dots, $a_{k-1}$ such that
$B_{k-1}=a_0b_0+a_1b_1+\dots+a_{k-1}b_{k-1}$,\quad
consequently we get
$pB_{k-1}=(a_0p)n+a_1(nc_1)+\dots+a_{k-1}(nc_{k-1})\equiv 0 \mod n$
and $p\equiv 0\mod\frac{n}{B_{k-1}}$ which proves the first part
of (5.1).

Suppose now that $p=\frac{n}{B_{k-1}}$. We have
$c_i = \frac{p}{n}b_i = \frac{b_i}{B_{k-1}}$
for ~$i=1$,\dots, $k-1$, hence
$\Gamma (f) =  \langle
  \frac{\be_1}{B_{k-1}},\dots,\frac{\be_{k-1}}{B_{k-1}}
  \rangle$.
\enddemo

Now, we can pass to the proof of (1.2).

Let $\phi=\phi(x,y)\in\series{x,y}$ be $y$--regular,
$(\phi,x)_0=\frac{n}{B_{k-1}}$. We shall check that
$(f,\phi)_0\leq\be_k$. If $k=h+1$ it is obvious
($\be_{h+1}={+}\infty$), so we assume $k\leq h$. Write
$$
  \phi=g_1\cdot\dots\cdot g_s,\qquad
   g_j\in\series{x,y} \text{ irreducible }
$$
We have
$$
\frac{(f,g_j)_0}{g_j,x)_0}\leq \frac{B_{k-1}\be_k}{n}\qquad
\text{for all } j=1,\dots,s                             \tag 1
$$
Indeed, if we had
$\frac{(f,g_j)_0}{(g_j,x)_0}>\frac{B_{k-1}\be_k}{n}$
for some $j$ then, by lemma ~5.1 we would get
$(g_j,x)_0\equiv 0\mod \frac{n}{B_k}$ and consequently
$(g_j,x)_0\geq \frac{n}{B_k}$.
It is impossible, because
$(g_j,x)_0\leq (\phi,x)_0=\frac{n}{B_{k-1}}$.

Now, from ~(1) we get
$$
(f,\phi)_0 =\sum_j(f,g_j)_0 \leq
 \sum_j\frac{B_{k-1}\be_k}{n}(g_j,x)_0 =
 \frac{B_{k-1}\be_k}{n}(\phi,x)_0 = \be_k
$$
Having proved the first part of (1.2) let us assume that
$(f,\phi)_0>n_{k-1}\be_{k-1}$. We claim that there exists a
$j\in\{\,1,\dots,s\,\}$ such that
$$
\frac{(f,g_j)_0}{(g_j,x)_0}>\frac{B_{k-2}\be_{k-1}}{n} \tag 2
$$
Suppose, contrary to our claim, that
$\frac{(f,g_j)_0}{(g_j,x)_0}\leq\frac{B_{k-2}\be_{k-1}}{n}$
for all $j=1$,\dots, $s$.
Thus we would have
$$
 (f,\phi)_0=\sum_j(f,g_j)_0 \leq
 \frac{B_{k-2}\be_{k-1}}{n}\sum_j(g_j,x)_0 =
 \frac{B_{k-2}\be_{k-1}}{n}(\phi,x)_0 = n_{k-1}\be_{k-1}
$$
which contradicts our assumption.

{}From (2) it follows, by lemma~5.1, that
$(g_j,x)_0=q\frac{n}{B_{k-1}}$ for some integer $q\geq 1$. On
the other hand $(g_j,x)_0\leq(\phi,x)_0=\frac{n}{B_{k-1}}$.
Therefore we get $q=1$ and $(g_j,x)_0=(\phi,x)_0$. Recall that
$g_j$ divides $\phi$, $g_j$ is irreducible and
$\ord g_j(0,y)=\ord \phi(0,y)$, thus $g_j$ is associated to
$\phi$ which proves irreducibility of $\phi$.


\enddocument